\documentclass[twoside,12pt]{article}
\usepackage[]{amsmath,amsthm,amssymb}
\usepackage[latin1]{inputenc}

\begin{document}

\title{{\Large\bf  A new index transform with the square of Whittaker's function}}

\author{Semyon  YAKUBOVICH}
\maketitle

\markboth{\rm \centerline{ Semyon   YAKUBOVICH}}{}
\markright{\rm \centerline{INDEX   TRANSFORM WITH THE SQUARE OF WHITTAKER'S FUNCTION}}

\begin{abstract} {\noindent An index transform, involving the square of Whittaker's function is introduced and investigated. The corresponding inversion formula is established. Particular cases cover index transforms of the Lebedev type with products of the modified Bessel functions.  }

\end{abstract}
\vspace{4mm}



\vspace{4mm}

\section {Introduction and preliminary results}

In 1964 Wimp [3] discovered a reciprocal pair of integral transformations, involving the Whittaker function $ W_{\mu,i\tau}(x)$ [1, Vol. III]

$$F(\tau)=  \int_0^\infty   W_{\mu,i\tau}(x)  f(x) {dx\over x^2},\  \tau >0,\eqno(1.1)$$

$$f(x)= { 1 \over \pi^2  } \int_0^\infty  \tau \sinh(2\pi\tau) \left|\Gamma\left({1\over 2}- \mu + i\tau\right)\right|^2  W_{\mu,i\tau}(x)  F(\tau) d\tau,\ x >0,\eqno(1.2)$$
 where $\mu \in \mathbb{R}$ and $\Gamma(z)$ is the Euler gamma function [1, Vol. III].    As we see in (1.1), (1.2), the corresponding integrals depend upon a second parameter (index) of the Whittaker function and involve it as the variable of integration.  Taking into account the known value
 
 $$W_{0,i\tau}(x) = \sqrt {{2x\over \pi}}\ K_{i\tau}(x),\eqno(1.3)$$ 
 where $K_\nu(x)$ is the modified Bessel or Macdonald function, we observe that the pair (1.1), (1.2) reduces to Kontorovich-Lebedev transforms [5].  Generally,  these functions are related through the following  integral (cf. Entry 2.16.8.4 in [1,  Vol. II])

 $$ \Gamma\left({1\over 2}- \mu + i\tau\right)\Gamma\left({1\over 2}- \mu - i\tau\right) W_{\mu,i\tau}(x) $$
 
 $$= 2(4x)^\mu e^{-x/2} \int_0^\infty  t^{-2\mu} e^{-t^2/(4x)} K_{2i\tau}(t) dt,\ x >0, \  {\rm Re} (\mu) < {1\over 2},\ \tau \in \mathbb{R}.\eqno(1.4)$$

This paper deals, however, with the index transformation, involving the square of the Whittaker function as the kernel

$$(F_\mu f)(x)=  {e^x\over x} \int_0^\infty   W^2_{\mu,i\tau}(x) f(\tau) d\tau,\quad x >0,\eqno(1.5)$$
where the real first index $\mu < 1/2$.  Our goal is to study its mapping properties and find an inversion formula.  The case $\mu=0$ corresponds to the Lebedev type index transform with the square of the Macdonald function

$$(F_0f)(x)=  { e^x \over \pi} \int_0^\infty  K^2_{i\tau}\left({x\over 2}\right) f(\tau) d\tau,\quad x >0.\eqno(1.6)$$
 The key ingredient for our departure point will be the integral representation of the square of Whittaker's function in terms of the Gauss hypergeometric function [1, Vol. III]
 (see [2,  formulas (3.26)]) 
 
 $$ W_{\mu, i\tau}^2 (x) = {x e^{-x}\over \Gamma(1-2\mu)}  \int_0^\infty t^{- 2\mu} e^{- x t} $$

$$\times \ {}_2F_1 \left( {1\over 2} -\mu-i\tau,\  {1\over 2} -\mu+ i\tau ; 1-2\mu ; - t^2- 2t \right) dt,\quad x >0,\ \mu < {1\over 2}.\eqno(1.7)$$
Meanwhile, the Gauss hypergeometric function can be represented by the integral (cf. [6]) 
 
 $${ \left|\Gamma(1/2-\mu +i\tau)\right|^2\over \Gamma(1-2\mu)} \  {}_2F_1 \left( {1\over 2} -\mu-i\tau,\  {1\over 2} -\mu+ i\tau ; 1-2\mu ;  - t^2 \right) $$

$$=  2 t^{2\mu} \int_0^\infty J_{-2\mu} (ty) K_{2i\tau} (y) dy.\eqno(1.8)$$
 Hence, employing the uniform estimate for the Macdonald function [5]
 
 $$\left|K_{i\tau}(x)\right| \le e^{-\delta |\tau|} K_0\left(x\cos\delta\right),\quad x >0, \tau \in \mathbb{R},\ \delta \in \left[0,  {\pi\over 2}\right),\eqno(1.9)$$ 
 and the known bound for the Bessel function $ |J_\nu(x)|\le C_\nu x^{-1/2} $ with some positive constant $C_\nu$, we get from (1.8) the inequality
 
 $${ \left|\Gamma(1/2-\mu +i\tau)\right|^2\over \Gamma(1-2\mu)} \  \left| {}_2F_1 \left( {1\over 2} -\mu-i\tau,\  {1\over 2} -\mu+ i\tau ; 1-2\mu ;  - t^2 \right) \right|$$

$$\le   2 C_{-2\mu}\   t^{2\mu-1/2}   e^{- 2\delta |\tau|}  \int_0^\infty y^{-1/2}  K_0\left(y\cos\delta\right)  dy$$

$$ =  {  C_{-2\mu}\  \Gamma^2(1/4) \over \sqrt 2 (\cos\delta)^{1/2}} \   t^{2\mu-1/2}   e^{- 2\delta |\tau|},\eqno(1.10)$$
 where the latter integral is calculated via relation 2.16.2.2 in [1, Vol. II]. Therefore for the square of the Whittaker function in (1.7) we obtain the following bound
 
 $$ W_{\mu, i\tau}^2 (x) \le  {  C_{-2\mu}\  \Gamma^2(1/4) \over \sqrt 2\  (\cos\delta)^{1/2}}\  {x e^{-x} e^{- 2\delta |\tau|} \over  \left|\Gamma(1/2-\mu +i\tau)\right|^2} $$
 
 $$\times   \int_0^\infty t^{- \mu-1/4} (t+2)^{\mu-1/4} e^{- x t}  dt \le  {  C_{-2\mu}\  \Gamma^2(1/4) \over \sqrt 2\  (\cos\delta)^{1/2}}$$
 
 $$\times   {x e^{-x} e^{- 2\delta |\tau|} \over  \left|\Gamma(1/2-\mu +i\tau)\right|^2}  \int_0^\infty t^{-1/2}  e^{- x t}  dt = {  C_{-2\mu}\  \Gamma^2(1/4) \over \sqrt 2\  (\cos\delta)^{1/2}}$$
 
 $$\times   {\sqrt{\pi x} \  e^{-x} e^{- 2\delta |\tau|} \over  \left|\Gamma(1/2-\mu +i\tau)\right|^2}, $$
i.e.

$$  W_{\mu, i\tau}^2 (x) \le   {  C_{-2\mu}\  \Gamma^2(1/4) \over \sqrt 2\  (\cos\delta)^{1/2}}  {\sqrt{\pi x} \  e^{-x} e^{- 2\delta |\tau|} \over  \left|\Gamma(1/2-\mu +i\tau)\right|^2},\eqno(1.11)$$
where $ x >0, \mu < 1/2, \tau \in \mathbb{R},\ \delta \in \left[0,  \pi/2\right)$.  This estimate allows to establish the composition structure of the index transform (1.5) in terms of the Laplace and Olevskii transforms [6]  under suitable condition on $f$.  In fact,  we have

{\bf Theorem 1}. {\it Let $x >0,\ \mu < 1/2, \ \delta \in \left[0,  \pi/2\right) $ and $f \in L_1\left( \mathbb{R}_+; e^{(\pi-2\delta) \tau} (\tau+1)^{2\mu} d\tau\right)$. Then the transform $(1.5)$ is a bounded map $F:  L_1\left( \mathbb{R}_+; \ e^{(\pi-2\delta) \tau} (\tau+1)^{2\mu} d\tau\right) \to C\left([x_0, \infty)\right),\\ x_0 > 0$ and can be represented as a composition of the Laplace and Olevskii transforms, namely,
 
$$(F_\mu f) (x)=    \int_0^\infty   t^{- 2\mu}  e^{- x t} G(t) dt,\eqno(1.12)$$
where}

$$G(t)= {1\over \Gamma(1-2\mu)}  \int_0^\infty  \ {}_2F_1 \left( {1\over 2} -\mu-i\tau,\  {1\over 2} -\mu+ i\tau ; 1-2\mu ; - t^2- 2t \right) f(\tau) d\tau.\eqno(1.13)$$
 
 \begin{proof} Indeed, the proof is immediate from the dominated convergence and Fubini's theorems which allow to interchange the order of integration after substitution integral (1.7) on the right-hand side of (1.5).  Moreover, since by virtue of (1.10)
 
 $$\left|(F_\mu f)(x)\right| \le  {1\over \Gamma(1-2\mu)} \int_0^\infty |f(\tau)|  \int_0^\infty t^{- 2\mu} e^{- x t} $$

$$\times \left|  {}_2F_1 \left( {1\over 2} -\mu-i\tau,\  {1\over 2} -\mu+ i\tau ; 1-2\mu ; - t^2- 2t \right) \right| dt d\tau $$

$$ \le   {  C_{-2\mu}\  \Gamma^2(1/4) \over \sqrt 2\  (\cos\delta)^{1/2}}  \int_0^\infty  { e^{- 2\delta \tau}\  |f(\tau)| \over  \left|\Gamma(1/2-\mu +i\tau)\right|^2}\ d\tau  \int_0^\infty  t^{- \mu-1/4} (t+2)^{\mu-1/4} e^{- x t}  dt $$

$$\le {  C_{-2\mu}\  \sqrt\pi\  \Gamma^2(1/4) \over \sqrt {2x_0} \  (\cos\delta)^{1/2}}  \int_0^\infty  { e^{- 2\delta \tau}\  |f(\tau)| \over  \left|\Gamma(1/2-\mu +i\tau)\right|^2}\ d\tau \eqno(1.14)$$
 and  due to the Stirling formula for the gamma function $\Gamma(1/2-\mu +i\tau) = O(\tau^{-\mu} e^{-\pi\tau/2} ),\ \tau \to + \infty$, we find that $F_\mu f$ is continuous on the interval $[x_0,\infty)$. Moreover, the norm estimate takes place 
 
 $$|| F_\mu f ||_{ C\left([x_0, \infty)\right)} = \sup_{ x \ge x_0} |F(x)| \le  {  C_{-2\mu}\  \sqrt\pi\  \Gamma^2(1/4) \over \sqrt {2x_0} \ (\cos\delta)^{1/2}}  \int_0^\infty  { e^{- 2\delta \tau}\  |f(\tau)| \over  \left|\Gamma(1/2-\mu +i\tau)\right|^2}\ d\tau $$
 
 $$ \le A_{\mu,\delta} || f||_{ L_1\left( \mathbb{R}_+;\  e^{(\pi-2\delta) \tau} (\tau+1)^{2\mu} d\tau\right)},$$
 where $A_{\mu,\delta} > 0$ is an absolute constant. 
 
 \end{proof}

 Writing the Olevskii transform (1.13) in the form
 
 $$G\left((t+1)^{1/2} -1\right)= {1\over \Gamma(1-2\mu)}  \int_0^\infty  \ {}_2F_1 \left( {1\over 2} -\mu-i\tau,\  {1\over 2} -\mu+ i\tau ; 1-2\mu ; - t \right) f(\tau) d\tau,\eqno(1.15)$$
 we appeal to the Plancherel theorem [6] for functions from the space $L_2(\mathbb{R}_+; |\Gamma(2i\tau)/ |\Gamma(1/2-\mu+ i\tau)|^2|^2 d\tau) $ which yields the Parseval equality
 
$$\int_0^\infty  \  { |f(\tau)|^2 \left|\Gamma(2i\tau)\right|^2 \over  |\Gamma(1/2-\mu+ i\tau)|^4} d\tau = {1\over 2\pi} \int_0^\infty \left| G\left((t+1)^{1/2} -1\right) \right|^2 t^{-2\mu} dt.\eqno(1.16)$$
 Moreover, the inverse operator is given by the formula
 $$ f(\tau)= \frac{ \left|\Gamma(1/2-\mu +i\tau)\right|^4} {2\pi \Gamma(1-2\mu) \left|\Gamma(2i\tau)\right|^2}  \int_0^\infty  t^{-2\mu}  {}_2F_1 \left( {1\over 2} -\mu-i\tau,\  {1\over 2} -\mu+ i\tau ; 1-2\mu ;  - t \right) $$
 
 $$\times G\left((t+1)^{1/2} -1\right) dt,\eqno(1.17)$$
where the convergence is with respect to the norm in $L_2(\mathbb{R}_+; |\Gamma(2i\tau)/ |\Gamma(1/2-\mu+ i\tau)|^2|^2 d\tau) $.  A simple change of variables permits to write (1.16), (1.17) in the equivalent form

$$\int_0^\infty  \  { |f(\tau)|^2 \left|\Gamma(2i\tau)\right|^2 \over  |\Gamma(1/2-\mu+ i\tau)|^4} d\tau = {1\over \pi} \int_0^\infty \left| G\left(t \right) \right|^2 (t+2)^{-2\mu} t^{-2\mu} (t+1) \ dt,\eqno(1.18)$$

 $$ f(\tau)= \frac{ \left|\Gamma(1/2-\mu +i\tau)\right|^4} {\pi \Gamma(1-2\mu) \left|\Gamma(2i\tau)\right|^2}  \int_0^\infty t^{-2\mu}  (t+2)^{-2\mu} \ (t+1)$$
 
 $$\times   {}_2F_1 \left( {1\over 2} -\mu-i\tau,\  {1\over 2} -\mu+ i\tau ; 1-2\mu ;  - t^2-2t \right) G\left( t\right) dt.\eqno(1.19)$$

\section{Inversion formula}

The main result of this section is the inversion formula for the transformation (1.5).   To establish this, we return to (1.12) and  apply the Mellin transform [4]  over a vertical straight line in the complex plane $s$ to both its sides under conditions of Theorem 1 to obtain

$$(F_\mu f)^*(s) \equiv \int_0^\infty (F_\mu f)(x) x^{s-1} dx =  \Gamma(s)   \int_0^\infty   t^{- 2\mu-s}   G(t) dt$$

$$ =   \Gamma(s)\  G^*(1-2\mu-s),\quad {\rm Re} (s) = \gamma \in \left({1\over 2},\   {3\over 4} -\ \mu\right),\ \mu < {1\over 4},\eqno(2.1)$$ 
 where the interchange of the order of integration is allowed due to the dominated convergence theorem. Indeed, from (1.13), (1.14) we observe
 
 $$|G(t)| \le  {1\over \Gamma(1-2\mu)} \int_0^\infty |f(\tau)| \left|  {}_2F_1 \left( {1\over 2} -\mu-i\tau,\  {1\over 2} -\mu+ i\tau ; 1-2\mu ; - t^2- 2t \right) \right| d\tau $$

$$ \le   {  C_{-2\mu}\  \Gamma^2(1/4)  \over \sqrt 2\  (\cos\delta)^{1/2}} \ t^{ \mu-1/4} (t+2)^{\mu-1/4}  \int_0^\infty  { e^{- 2\delta \tau}\  |f(\tau)| \over  \left|\Gamma(1/2-\mu +i\tau)\right|^2}\ d\tau \eqno(2.2)$$
and the integral 

$$\int_0^\infty   t^{- \mu-s -1/4} (t+2)^{\mu-1/4}   dt$$
converges  when $3/4-\mu > \gamma > 1/2$.   In the mean time, reciprocally via the inverse Mellin transform in $L_2$ we have

$$ t^{2\mu-1} G \left({1\over t}\right) = {1\over 2\pi i} \int_{\gamma -i\infty}^{\gamma+i\infty} { (F_\mu f)^*(s)\over \Gamma(s)}\  t^{-s} ds,\quad t > 0\eqno(2.3)$$
which is equivalent to the equality

 $$ G \left( t\right) = {1\over 2\pi i} \int_{1-2\mu-\gamma -i\infty}^{1-2\mu-\gamma+i\infty} { (F_\mu f)^*(1-2\mu-s)\over \Gamma(1-2\mu-s)}\  t^{-s} ds.\eqno(2.4)$$
But the Parseval equality for the Mellin transform in $L_2$ [4] implies

$$\int_{0}^\infty |G(t)|^2 t^{-4\mu-2\gamma +1} dt = {1\over 2\pi}  \int_{-\infty}^{\infty} \left|   { F^*(\gamma+iu)\over \Gamma(\gamma+ iu)} \right|^2 du.\eqno(2.5)$$
In fact, when  $ \gamma \in \left(1/2,\   3/4 -\ \mu\right)$ the integral on the left-hand side of (2.5) converges.   Precisely, it has via (2.2)

$$\int_{0}^\infty |G(t)|^2 t^{-4\mu-2\gamma +1} dt  = \left(\int_0^1 + \int_1^\infty \right)  |G(t)|^2 t^{-4\mu-2\gamma +1} dt  $$

$$\le A \int_0^1 t^{-2\mu-2\gamma + 1/2} dt + B  \int_1^\infty  t^{-2\gamma} dt < \infty, $$
 where $A, B$ are positive constants.  Therefore, the convergence of the integral (2.4) in the mean square sense with respect to the norm in  $L_2(\mathbb{R}_+;  t^{1-4\mu-2\gamma } dt )$. Moreover,  assuming the condition 
 
 $$f \in L_1\left( \mathbb{R}_+; e^{(\pi-2\delta) \tau} (\tau+1)^{2\mu} d\tau\right) \cap L_2\left(\mathbb{R}_+; {|\Gamma(2i\tau)|^2 d\tau \over  |\Gamma(1/2-\mu+ i\tau)|^4 } \right),\eqno(2.6)$$
 where $\delta \in [0, \pi/2)$, it yields immediately (see (1.18), (2.2)) that
 
$$G\in   L_2\left( \mathbb{R}_+;   (t+2)^{-2\mu} t^{-2\mu} (t+1) d t\right) \cap L_2\left( \mathbb{R}_+;  t^{-4\mu-2\gamma +1} dt \right),\eqno(2.7)$$
 where $\mu < 1/4, \  \gamma \in \left(1/2,\   3/4 -\ \mu\right)$. Therefore, returning to (1.19), we  write it as the $L_2$-limit $f= \lim_{N\to \infty} f_N$ with respect to the norm in $L_2(\mathbb{R}_+; |\Gamma(2i\tau)/ |\Gamma(1/2-\mu+ i\tau)|^2|^2 d\tau)$, where
 
 $$ f_N(\tau)= \frac{ \left|\Gamma(1/2-\mu +i\tau)\right|^4} {\pi \Gamma(1-2\mu) \left|\Gamma(2i\tau)\right|^2}  \int_{1/N}^N  t^{-2\mu} (t+2)^{-2\mu} (t+1)$$

$$\times   {}_2F_1 \left( {1\over 2} -\mu-i\tau,\  {1\over 2} -\mu+ i\tau ; 1-2\mu ;  - t^2-2t  \right) G(t) dt.\eqno(2.8)$$ 
In the meantime, Entry 8.4.49.24 in [1, Vol. III]  gives the following Mellin-Barnes integral representation for the Gauss hypergeometric function

$$ (2t+t^2)^{-2\mu}  {}_2F_1 \left( {1\over 2} -\mu-i\tau,\  {1\over 2} -\mu+ i\tau ; 1-2\mu ;  - 2t- t^2 \right)  $$

$$= {\Gamma(1-2\mu) \over 2\pi i }  \int_{\alpha-i\infty}^{\alpha+i\infty}  \Gamma\left( {1\over 2} +\mu-i\tau-w \right) \Gamma\left( {1\over 2} +\mu+i\tau -w\right)  { (t+1)^{-2w}\over \Gamma ^2( 1 - w) }  dw,\eqno(2.9)$$
where $t >0, \ \alpha < 1/2+\mu.$  Meanwhile, it is possible to replace in (2.9)  the vertical line $(\alpha-i\infty,\ \alpha+i\infty)$ in the complex plane $w$ on the right-hand loop $\cal{L}_{+\infty}$ with a  positive orientation and bounded imaginary part which comprises   right-hand simple poles of gamma functions $w= 1/2+\mu \pm i\tau + n,\ n \in \mathbb{N}$, separating residues at two simple poles $w= 1/2+\mu \pm i\tau$. In fact, it is possible via analyticity property of the integrand and its exponential decay when $t > 0$ and ${\rm Re} (w) \to +\infty$ owing to the Stirling asymptotic formula for the gamma function.  Hence plugging the obtained  integral in (2.8), it gives 

 $$ f_N(\tau)= \frac{ \left [\Gamma(1/2-\mu -i\tau)\right]^2} {\pi\ \Gamma(-2i\tau)}  \int_{1/N}^N  {G(t) \ dt \over (t+1)^{2(\mu-i\tau)} } $$

$$+ \frac{ \left [\Gamma(1/2-\mu+i\tau)\right]^2} {\pi\ \Gamma(2i\tau)}  \int_{1/N}^N  {G(t) \ dt \over (t+1)^{2(\mu+i\tau)} } $$

$$ - \frac{ \left|\Gamma(1/2-\mu +i\tau)\right|^4} {2\pi^2 i \left|\Gamma(2i\tau)\right|^2}  \int_{1/N}^N   G(t) $$
 
 $$\times \int_{\cal{L}_{+}}  \Gamma\left( {1\over 2} +\mu-i\tau-w \right) \Gamma\left( {1\over 2} +\mu+i\tau -w\right)  {(t+1)^{1-2w} \over  \Gamma ^2( 1 - w) } \ dw dt.\eqno(2.10)$$ 
But since
$$ (t+1)^{1-2w} = {1\over \Gamma(2w-1)} \int_0^\infty e^{-y(t+1)} y^{2(w-1)} dy,\ {\rm Re} (w) > {1\over 2},\eqno(2.11)$$
we  have 

$$   f_N(\tau)=  \frac{ \left [\Gamma(1/2-\mu -i\tau)\right]^2} {\pi\ \Gamma(-2i\tau)}  \int_{1/N}^N  {G(t) \ dt \over (t+1)^{2(\mu-i\tau)} } $$

$$+ \frac{ \left [\Gamma(1/2-\mu+i\tau)\right]^2} {\pi\ \Gamma(2i\tau)}  \int_{1/N}^N  {G(t) \ dt \over (t+1)^{2(\mu+i\tau)} }+ I_N(\tau), \eqno(2.12)$$
where

$$I_N(\tau)= - \frac{ \left|\Gamma(1/2-\mu +i\tau)\right|^4  } {2\pi^2 i\ \left|\Gamma(2i\tau)\right|^2}  \int_{1/N}^N   \int_{\cal{L}_{+}}  {\Gamma\left( {1\over 2} +\mu-i\tau-w \right) \Gamma\left( {1\over 2} +\mu+i\tau -w\right)\over  \Gamma ^2( 1 - w) \ \Gamma(2w-1)}  $$

$$\times  G(t) \int_0^\infty e^{-y(t+1)} y^{2(w-1)}\  dy dw dt.\eqno(2.13)$$ 
So, we are ready to interchange the order of integration in the latter iterated integral via the dominated convergence theorem under the condition (2.7). Consequently,

$$I_N(\tau) =- \frac{ \left|\Gamma(1/2-\mu +i\tau)\right|^4  } {2\pi^2 i\ \left|\Gamma(2i\tau)\right|^2}    \int_{\cal{L}_{+\infty}}  {\Gamma\left( {1\over 2} +\mu-i\tau-w \right) \Gamma\left( {1\over 2} +\mu+i\tau -w\right)\over  \Gamma ^2( 1 - w) \ \Gamma(2w-1)}  $$

$$\times \int_0^\infty e^{-y} y^{2(w-1)}\  dy  \int_{1/N}^N e^{-yt}  G(t)  dt dy  dw.\eqno(2.14)$$ 
Moreover, it is possible to pass to the limit under the integral sign in (2.14) since the Cauchy-Schwarz inequality  implies

$$\int_{1/N}^N e^{-yt}\  |G(t)| \  dt \le ||G||_{L_2\left( \mathbb{R}_+;   (t+2)^{-2\mu} t^{-2\mu} (t+1) dt \right)} \left(\int_0^\infty  e^{-2yt} \  {(t+2)^{2\mu}\  t^{2\mu} \over t+1}\ dt  \right)^{1/2} $$

$$ \le  2^\mu ||G||_{L_2\left( \mathbb{R}_+;   (t+2)^{-2\mu} t^{-2\mu} (t+1) dt \right)} \left(\int_0^\infty  e^{-2yt} \  (t+1)^{2\mu-1}\  t^{2\mu}  dt  \right)^{1/2} $$

$$= 2^\mu ||G||_{L_2\left( \mathbb{R}_+;   (t+2)^{-2\mu} t^{-2\mu} (t+1) dt \right)}\   \Psi^{1/2} (2\mu+1, 4\mu+1; 2y),\quad -{1\over 2} < \mu < {1\over 4},$$
where $\Psi (2\mu+1,\ 4\mu+1; 2y )$ is the Tricomi function which behaves as $O(y^{-4\mu}), \ \mu\neq 0, \  O(\log y),\ \mu=0, \ y \to 0$ and $O(y^{-2\mu-1}),\ y \to \infty$. It yields that  the integral

$$\int_0^\infty e^{-y} y^{2 (w - 1)}\   \Psi^{1/2} (2\mu+1, 4\mu+1; 2y) dy$$
converges when $ {\rm Re} (w) > {1\over 2} + \mu.$  Thus we derive 

$$I(\tau) = \lim_{N\to \infty} I_N(\tau)=  - \frac{ \left|\Gamma(1/2-\mu +i\tau)\right|^4  } {2\pi^2 i\ \left|\Gamma(2i\tau)\right|^2}   \int_{\cal{L}_{+\infty}}  {\Gamma\left( {1\over 2} +\mu-i\tau-w \right) \Gamma\left( {1\over 2} +\mu+i\tau -w\right)\over  \Gamma ^2( 1 - w) \ \Gamma(2w-1)}  $$

$$\times  \int_0^\infty e^{-y} y^{2(w-1)}\    \int_{0}^\infty e^{-yt}  G(t)  dt dy dw.\eqno(2.15) $$
But via (2.7) we find, in turn,  

$$\int_0^\infty e^{-y} y^{2( {\rm Re} (w) -1)}\    \int_{0}^\infty e^{-yt}  |G(t)|  dt dy \le  ||G||_{L_2\left( \mathbb{R}_+; \  t^{1-4\mu-2\gamma} dt \right)}  \int_0^\infty e^{-y} y^{2( {\rm Re} (w) -1)}\  $$

$$\times \left(  \int_{0}^\infty e^{-2yt}  t^{4\mu+2\gamma-1} dt\right)^{1/2} dy =  2^{-2\mu-\gamma}\ \Gamma^{1/2} (4\mu+2\gamma)\ ||G||_{L_2\left( \mathbb{R}_+; \  t^{1-4\mu-2\gamma} dt \right)} $$

$$\times   \int_0^\infty e^{-y} y^{2( {\rm Re} (w) -1-\mu)-\gamma} dy =  2^{-2\mu-\gamma}\ \Gamma^{1/2} (4\mu+2\gamma)\  \Gamma(2( {\rm Re} (w)-\mu)-\gamma-1) $$

$$\times  ||G||_{L_2\left( \mathbb{R}_+; \  t^{1-4\mu-2\gamma} dt \right)} < \infty $$
when $ {\rm Re} (w)  > (\gamma+1)/2 +\mu$.   Therefore, returning to (2.15) and taking into account (2.4), we appeal to the Mellin-Parseval equality [4] to obtain

$$I(\tau) = - \frac{ \left|\Gamma(1/2-\mu +i\tau)\right|^4  } {(2\pi i)^2\ \pi  \left|\Gamma(2i\tau)\right|^2}   \int_{\cal{L}_{+\infty}}  {\Gamma\left( {1\over 2} +\mu-i\tau-w \right) \Gamma\left( {1\over 2} +\mu+i\tau -w\right)\over  \Gamma ^2( 1 - w) \ \Gamma(2w-1)}  $$

$$\times  \int_0^\infty e^{-y} y^{2(w-1)}\   \int_{1-2\mu-\gamma -i\infty}^{1-2\mu-\gamma+i\infty} { (F_\mu f)^*(1-2\mu-s) \ \Gamma (1-s) \over \Gamma(1-2\mu-s)}\  y^{s-1} ds dy dw$$

 $$= - \frac{ \left|\Gamma(1/2-\mu +i\tau)\right|^4  } {(2\pi i)^2\ \pi  \left|\Gamma(2i\tau)\right|^2}   \int_{\cal{L}_{+\infty}}  {\Gamma\left( {1\over 2} +\mu-i\tau-w \right) \Gamma\left( {1\over 2} +\mu+i\tau -w\right)\over  \Gamma ^2( 1 - w) \ \Gamma(2w-1)}  $$

$$\times   \int_{1-2\mu-\gamma -i\infty}^{1-2\mu-\gamma+i\infty} { (F_\mu f)^*(1-2\mu-s) \ \Gamma (1-s) \over \Gamma(1-2\mu-s)}\  \Gamma (2w +s-2) ds  dw.\eqno(2.16)$$
Moreover, the interchange of the order of integration is allowed on the latter iterated integral in (2.15) because the imaginary part  of $w$ is bounded, and we deduce

$$I(\tau) = - \frac{ \left|\Gamma(1/2-\mu +i\tau)\right|^4  } {2\pi^2 i  \left|\Gamma(2i\tau)\right|^2} \int_{1-2\mu-\gamma -i\infty}^{1-2\mu-\gamma+i\infty} { (F_\mu f)^*(1-2\mu-s) \ \Gamma (1-s) \over \Gamma(1-2\mu-s)}\  \Phi_{\mu,i\tau}(s) ds,\eqno(2.17)$$  
where 

$$\Phi_{\mu,i\tau} (s) = {1\over 2\pi i}  \int_{\cal{L}_{+\infty}}  {\Gamma\left( {1\over 2} +\mu-i\tau-w \right) \Gamma\left( {1\over 2} +\mu+i\tau -w\right) \Gamma (2w +s-2) \over  \Gamma ^2( 1 - w) \ \Gamma(2w-1)} dw.\eqno(2.18)$$
Further, involving the beta integral, we write (2.17) as follows

$$I(\tau) = - \frac{ \left|\Gamma(1/2-\mu +i\tau)\right|^4  } {2\pi^2 i\   \Gamma(1-2\mu) \left|\Gamma(2i\tau)\right|^2} \int_{1-2\mu-\gamma -i\infty}^{1-2\mu-\gamma+i\infty}  (1-2\mu-s) (F_\mu f)^*(1-2\mu-s)\  \Phi_{\mu, i\tau} (s) $$

$$\times \int_0^1 x^{-s} (1-x)^{-2\mu} dx ds= - \frac{ \left|\Gamma(1/2-\mu +i\tau)\right|^4  } {2\pi^2 i\  \Gamma(1-2\mu) \left|\Gamma(2i\tau)\right|^2} \int_0^1 (1-x)^{-2\mu} $$

$$\times \int_{1-2\mu-\gamma -i\infty}^{1-2\mu-\gamma+i\infty}  (1-2\mu-s) (F_\mu f)^*(1-2\mu-s)\  \Phi_{\mu,i\tau} (s)  x^{-s} ds dx$$  

$$= - \frac{ \left|\Gamma(1/2-\mu +i\tau)\right|^4  } {2\pi^2 i  \Gamma(1-2\mu) \left|\Gamma(2i\tau)\right|^2} \int_0^1 (1-x)^{-2\mu} x^{2\mu} $$

$$\times {d\over dx} \int_{1-2\mu-\gamma -i\infty}^{1-2\mu-\gamma+i\infty} (F_\mu f)^*(1-2\mu-s)\  \Phi_{\mu, i\tau} (s)  x^{1-2\mu-s} ds dx,\eqno(2.19)$$ 
where the interchange of the order of integration and the differentiation under the integral sign are permitted by the absolute and uniform convergence.  Indeed, we have the estimate for the corresponding iterated integral (see (2.4), (2.18))

$$\int_{1-2\mu-\gamma -i\infty}^{1-2\mu-\gamma+i\infty}  \left| (1-2\mu-s) (F_\mu f)^*(1-2\mu-s)\  \Phi_{\mu, i\tau} (s) \right|$$

$$\times \int_0^1 x^{2\mu+\gamma-1} (1-x)^{-2\mu} dx |ds| \le B(2\mu+\gamma, 1-2\mu) \left(\int_{1-2\mu-\gamma -i\infty}^{1-2\mu-\gamma+i\infty} \left| { (F_\mu f)^*(1-2\mu-s)\over \Gamma(1-2\mu-s)} \right|^2 \right)^{1/2} $$

$$\times \left(\int_{1-2\mu-\gamma -i\infty}^{1-2\mu-\gamma+i\infty} \bigg| \Gamma (2-2\mu-s)\  \Phi_{\mu, i\tau} (s) \bigg|^2 \right)^{1/2} < \infty,$$
where $B(a,b)$ is the Euler beta function.  Furthermore, it yields that $(F_\mu f)^*(1-2\mu-s) \in L_2(1-2\mu-\gamma-i\infty,\  1-2\mu-\gamma+i\infty)$ and the Mellin-Parseval equality [4] presumes the representation from (2.18)

$$I(\tau) = - \frac{ \left|\Gamma(1/2-\mu +i\tau)\right|^4  } {\pi \Gamma(1-2\mu) \left|\Gamma(2i\tau)\right|^2} \int_0^1 (1-x)^{-2\mu} x^{2\mu}$$

$$\times  {d\over dx} \bigg[ x^{1-2\mu} \int_0^\infty y^{-2\mu} \ (F_\mu f) (y)\  \varphi_{\mu,i\tau} (xy) dy \bigg] dx,\eqno(2.20)$$ 
where  $(F_\mu f) (y) \in L_2 \left(\mathbb{R}_+;  y^{2\gamma -1} dy\right)$ and 
$$\varphi_{\mu,i\tau} (x) = {1\over 2\pi i}  \int_{1-2\mu-\gamma -i\infty}^{1-2\mu-\gamma+i\infty}  \Phi_{\mu,i\tau} (s)  x^{-s} ds.\eqno(2.21)$$
Returning to (2.18), we substitute its right-hand side in  (2.21) and interchange the order of integration by Fubini's theorem.  Hence after calculation the inner Mellin-Barnes integral it immediately implies
$$ \varphi_{\mu,i\tau} (x) = {e^{-x} \over 2\pi i}   \int_{\cal{L}_{+\infty}}  {\Gamma\left( {1\over 2} +\mu-i\tau-w \right) \Gamma\left( {1\over 2} +\mu+i\tau -w\right) \over  \Gamma ^2( 1 - w) \ \Gamma(2w-1)}\  x^{2(w-1)} dw.\eqno(2.22)$$
Employing the residue theorem,  we  calculate the kernel $\varphi_{\mu, i\tau}$ in terms of the generalized hypergeometric functions.  Indeed, we derive

$$\varphi_{\mu, i\tau} (x) = e^{-x} x^{2\mu-1} \bigg[ \sum_{n=1}^\infty  \frac{ (-1)^{n+1} x^{2(n-i\tau)}\  \Gamma(2i\tau-n) }{ n!\ \Gamma^2( 1/2-\mu +i\tau-n) \Gamma(2(\mu-i\tau+n))}\bigg. $$

$$\bigg. + \sum_{n=1}^\infty  \frac{ (-1)^{n+1} x^{2(n+i\tau)}\  \Gamma(-2i\tau-n) }{ n!\  \Gamma^2( 1/2-\mu -i\tau-n) \Gamma(2(\mu+i\tau+n))}  \bigg]$$

$$= - {e^{-x} (x/2)^{2\mu-1}\over \sqrt \pi \sinh(2\pi\tau)}  \bigg[ \cos^2(\pi(\mu-i\tau)) \left({2\over x} \right)^{2i\tau} \sum_{n=1}^\infty  \frac{  \Gamma ( 1/2+\mu -i\tau+n) }{ n!\   \Gamma(1-2i\tau+n)  \Gamma(\mu-i\tau+n)} \left({x^2\over 4}\right)^n \bigg. $$

$$\bigg. - \cos^2(\pi(\mu+i\tau)) \left({2\over x} \right)^{-2i\tau} \sum_{n=1}^\infty  \frac{ \Gamma( 1/2+\mu +i\tau+n) }{ n!\  \Gamma(1+2i\tau+n)  \Gamma(\mu+i\tau+n)} \left({x^2\over 4}\right)^n \bigg]$$

$$= -   e^{-x} x^{2\mu-1} \bigg[  \frac{ x^{-2i\tau}\   \Gamma(2i\tau)  }{  \Gamma^2 ( 1/2-\mu +i\tau)\  \Gamma(2(\mu-i\tau))} $$

$$\times \left[ {}_1F_2 \left( {1\over 2} +\mu -i\tau ;  \ 1-2i\tau, \   \mu-i\tau;\  {x^2\over 4} \right) -1\right]  $$

$$\bigg. + \frac{ x^{2i\tau}\ \Gamma(-2i\tau)  }{  \Gamma^2 ( 1/2-\mu -i\tau)\  \Gamma(2(\mu+i\tau))}$$

$$\times \left[ {}_1F_2 \left( {1\over 2} +\mu +i\tau ;  \ 1+2i\tau, \   \mu+i\tau;\  {x^2\over 4} \right) -1\right] \bigg]. $$ 
Hence we establish the following equality

$$ \varphi_{\mu, i\tau} (x) = -   2 e^{-x} x^{2\mu-1}\    {\rm Re} \bigg[  \frac{  x^{-2i\tau}\  \Gamma(2i\tau)  }{  \Gamma^2 ( 1/2-\mu +i\tau)\  \Gamma(2(\mu-i\tau))} \bigg.$$

$$\times \left[ {}_1F_2 \left( {1\over 2} +\mu -i\tau ;  \ 1-2i\tau, \   \mu-i\tau;\  {x^2\over 4} \right) -1\right] \bigg],\quad x >0.\eqno(2.23) $$ 
Further, fulfilling the differentiation  in (2.20) owing to the absolute and uniform convergence for each  $x >0$, we find

$$I(\tau) = - \frac{ \left|\Gamma(1/2-\mu +i\tau)\right|^4  } {\pi \Gamma(1-2\mu) \left|\Gamma(2i\tau)\right|^2} \bigg[ (1-2\mu) \int_0^1 (1-x)^{-2\mu}  \int_0^\infty y^{-2\mu} \ (F_\mu f) (y)\  \varphi_{\mu,i\tau} (xy) dy dx \bigg.$$

$$\bigg. +  \int_0^1 x  (1-x)^{-2\mu}  \int_0^\infty y^{1-2\mu} \ (F_\mu f) (y)\  \varphi_{\mu,i\tau}^\prime (xy) dy dx\bigg].\eqno(2.24)$$ 
Moreover, we use the integration  by parts  in the integral by $y$ of the first term on the right-hand side of (2.24), eliminating the  integrated terms with the aim of the estimate (1.14) and since (cf. [1, Vol. III])

$$ \varphi_{\mu,i\tau} (y) = O( y^{2\mu +1}),\ y \to 0,\quad  \varphi_{\mu,i\tau} (y) = O( y^{2\mu -1}),\ y \to \infty.$$
Consequently, equality (2.24) becomes

$$I(\tau) =  \frac{ \left|\Gamma(1/2-\mu +i\tau)\right|^4  } {\pi \Gamma(1-2\mu) \left|\Gamma(2i\tau)\right|^2}  \int_0^1 (1-x)^{-2\mu}  \int_0^\infty y^{1-2\mu} \  \varphi_{\mu,i\tau} (xy) \ (F_\mu f)^\prime  (y) dy dx. \eqno(2.25)$$
The interchange of the order of integration in (2.5) allows to write it in the form, employing the Riemann-Liouville fractional integral of the order $1-2\mu$  [4]

$$I(\tau) =  \frac{ \left|\Gamma(1/2-\mu +i\tau)\right|^4  } {\pi  \left|\Gamma(2i\tau)\right|^2}   \int_0^\infty  I_{0+}^{1-2\mu} \left( \varphi_{\mu,i\tau} \right)(y) \ (F_\mu f)^\prime  (y) dy dx, \eqno(2.26)$$
where 
$$ \left(I_{0+}^{1-2\mu}  \varphi_{\mu, i\tau} \right) (y) =  {1\over \Gamma(1-2\mu)}  \int_0^y (y-x)^{-2\mu} \varphi_{\mu, i\tau} (x) dx.\eqno(2.27)$$
In fact, the interchange is guaranteed by the estimate  (see (1.12))

$$\int_0^1 (1-x)^{-2\mu}  \int_0^\infty y^{1-2\mu} \  \left|\varphi_{\mu,i\tau} (xy) \ (F_\mu f)^\prime  (y)\right| dy dx$$

$$\le \int_0^1 (1-x)^{-2\mu}  \int_0^\infty y^{1-2\mu} \  \left|\varphi_{\mu,i\tau} (xy)\right|   \int_0^\infty   t^{1- 2\mu}  e^{- y t}  |G(t)| \ dt dy dx$$

$$\le  ||G||_{L_2\left( \mathbb{R}_+; \  t^{1-4\mu-2\gamma} dt \right)} \int_0^1 (1-x)^{-2\mu}  \int_0^\infty y^{1-2\mu} \  \left|\varphi_{\mu,i\tau} (xy)\right|   \left (\int_0^\infty   t^{1+2\gamma}  e^{- 2y t}  dt \right)^{1/2} dy dx$$

$$=  ||G||_{L_2\left( \mathbb{R}_+; \  t^{1-4\mu-2\gamma} dt \right)} \int_0^1 (1-x)^{-2\mu}  x^{2\mu+\gamma-1} dx \int_0^\infty y^{-2\mu-\gamma} \  \left|\varphi_{\mu,i\tau} (y)\right|  dy $$

$$\times \left (\int_0^\infty   t^{1+2\gamma}  e^{- 2 t}  dt \right)^{1/2}  < \infty$$
when $2\mu+\gamma >0$ by virtue of the dominated convergence theorem.  It holds, when (see above )  $ \hbox{max} (-2\mu, 1/2) < \gamma <  3/4-\mu,\ \mu \in (-1/2,\ 1/4).$

We summarize these results by the following theorem.

{\bf Theorem 2}. {\it Let $\tau >0,\ -1/2 < \mu < 1/4$ and $f(\tau)$ satisfy condition $(2.6)$ which implies $(2.7)$ with $\gamma \in (\hbox{max}\ (-2\mu, 1/2),\  3/4-\mu)$. Then the sequence $\{I_N(\tau)\}_{N \ge 1}$,  being defined in $(2.13)$,  converges pointwisely to $I(\tau)$ for all $\tau > 0$ and can be represented by the absolutely convergent integral $(2.26)$, where $F_\mu f$ is defined by $(1.5)$ and the kernel $\varphi_{\mu, i\tau}$ by the formula  $(2.23)$}.

Now we recall (2.12) and pass to the limit when $N \to \infty$ with respect to the norm in $L_2(\mathbb{R}_+; |\Gamma(2i\tau)/ |\Gamma(1/2-\mu+ i\tau)|^2|^2 d\tau) $.  But under conditions of Theorem 2 the sequence $\{I_N(\tau)\}_{N \ge 1}$ converges to the same limit  $I(\tau)$ given by (2.26) pointwisely.  Therefore it gives

$$   f(\tau)=  {1\over \pi}\   \frac{ \Gamma^2 (1/2-\mu -i\tau)} { \Gamma(-2i\tau)}  \int_{0}^\infty  {G(t) \ dt \over (t+1)^{2(\mu-i\tau)} } $$

$$+  {1\over \pi}\   \frac{ \Gamma^2 (1/2-\mu +i\tau)} { \Gamma(2i\tau)}  \int_{0}^\infty  {G(t) \ dt \over (t+1)^{2(\mu+i\tau)} } $$

$$ + \frac{ \left|\Gamma(1/2-\mu +i\tau)\right|^4  } {\pi  \left|\Gamma(2i\tau)\right|^2}\int_0^\infty  \left( F_\mu f\right)^\prime(y) \left(I_{0+}^{1-2\mu} \varphi_{\mu, i\tau}\right) (y) dy.\eqno(2.28)$$ 
 The first two integrals on the right-hand side of (2.28) converges with respect to the norm in $L_2(\mathbb{R}_+; |\Gamma(2i\tau)/ |\Gamma(1/2-\mu+ i\tau)|^2|^2 d\tau) $. Indeed, appealing to (1.18) and Parseval's equality for the Fourier transform, we denote by  
 
 $$g_N(\tau)= \frac{ \Gamma^2 (1/2-\mu -i\tau)} {\pi\ \Gamma(-2i\tau)}  \int_{1/N}^N  {G(t) \ dt \over (t+1)^{2(\mu-i\tau)} }\eqno(2.29)$$
and its limit by $g(\tau)$  to derive
 
 $$||g_N-g||_{L_2(\mathbb{R}_+; |\Gamma(2i\tau)/ |\Gamma(1/2-\mu+ i\tau)|^2|^2 d\tau) } = {1\over \pi} \left(\int_0^\infty    \bigg| \left( \int_0^{1/N} + \int_N^\infty \right)  {G(t) \ dt \over (t+1)^{2(\mu-i\tau)} } \bigg|^2 d\tau\right)^{1/2}$$
 
 $$= {1\over 2\pi} \left(\int_0^\infty    \bigg| \left( \int_0^{2\log(1/N+1)} + \int_{2\log N}^\infty \right)  e^{(1/2-\mu) u} \ G\left(e^{u/2} -1\right) e^{i\tau u} \ du \bigg|^2 d\tau\right)^{1/2}$$

 $$\le C  \left( \int_0^{2\log(1/N+1)}  e^{(1-2\mu) u} \ \left|G\left(e^{u/2} -1\right)\right|^2 \ du \right)^{1/2} $$
 
 $$+ C \left( \int_{2\log N}^\infty  e^{(1-2\mu) u} \ \left|G\left(e^{u/2} -1\right)\right|^2 \ du \right)^{1/2}$$ 
 
$$= \sqrt 2  C  \left( \int_0^{1/N}  (t+1)^{1-4\mu} \ \left|G\left(t\right)\right|^2 \ dt \right)^{1/2} $$

$$+ \sqrt 2 C \left( \int_{ N}^\infty   (t+1)^{1-4\mu} \ \left|G\left(t\right)\right|^2 \ dt  \right)^{1/2}\  \to 0,\ N \to \infty,$$ 
 where $C >0$ is a constant. Analogously the second integral on the right-hand side in (2.28) can be treated.  On the other hand,   taking  into account (2.14), (2.22),  (2.29), we  write 
 
 $$   f_N(\tau)=  g_N(\tau) + g_N(-\tau)  - \frac{ \left|\Gamma(1/2-\mu +i\tau)\right|^4  } {\pi\ \left|\Gamma(2i\tau)\right|^2 \Gamma(1-2\mu)} $$
 
 $$\times  \int_0^\infty  \varphi_{\mu, i\tau} (y) \int_0^\infty x^{-2\mu} \int_{1/N}^N t^{1-2\mu}  e^{-(x+y)t}  G(t)\  dt dx dy$$ 

$$ =  g_N(\tau) + g_N(-\tau)  - \frac{ \left|\Gamma(1/2-\mu +i\tau)\right|^4  } {\pi\ \left|\Gamma(2i\tau)\right|^2 } $$
 
 $$\times  \int_0^\infty  \varphi_{\mu, i\tau} (y) \left(D_{-}^{2\mu} F_{\mu, N} f \right) (y)  dy,$$ 
where

$$ \left(D_{-}^{2\mu}  f\right)(x) = - {1\over \Gamma(1-2\mu)} \int_x^\infty (y-x)^{-2\mu}  f^\prime(y) dy $$
is the right-sided fractional derivative of the order $2\mu$ and  (see (1.12))

 $$(F_{\mu, N} f)(x)=  \int_{1/N}^N t^{-2\mu}  e^{-xt}  G(t)  dt.$$
 Then plugging in the value (2.23) of the kernel $\varphi_{\mu, i\tau}$, we have
 
 $$   f_N(\tau)= g_N(\tau) + g_N(-\tau)  -   \frac{  \Gamma^2 ( 1/2-\mu -i\tau) }{ \pi  \Gamma(-2i\tau)\  \Gamma(2(\mu-i\tau))}  $$
 
 $$\times  \int_0^\infty     e^{-y} y^{2(\mu-i\tau)-1}  \left(D_{-}^{2\mu} F_{\mu, N} f \right) (y)  dy$$ 
 
$$ -   \frac{  \Gamma^2 ( 1/2-\mu +i\tau) }{ \pi  \Gamma(2i\tau)\  \Gamma(2(\mu+i\tau))}    \int_0^\infty     e^{-y} y^{2(\mu+i\tau)-1}  \left(D_{-}^{2\mu} F_{\mu, N} f \right) (y)  dy$$ 
 
 $$+ {2\over \pi} \int_0^\infty      {\rm Re} \bigg[  \frac{  \Gamma^2 ( 1/2-\mu -i\tau)    }{  \Gamma(-2i\tau)  \Gamma(2(\mu-i\tau))}  {}_1F_2 \left( {1\over 2} +\mu -i\tau ;  \ 1-2i\tau, \   \mu-i\tau;\  {y^2\over 4} \right) \ y^{2(\mu-i\tau)-1} \bigg] $$

$$\times e^{-y} \left(D_{-}^{2\mu} F_{\mu, N} f \right) (y)  dy.$$
In the meantime, we observe 

$$\frac{  \Gamma^2 ( 1/2-\mu -i\tau) }{ \pi  \Gamma(-2i\tau)\  \Gamma(2(\mu-i\tau))}    \int_0^\infty     e^{-y} y^{2(\mu-i\tau)-1}  \left(D_{-}^{2\mu} F_{\mu, N} f \right) (y)  dy$$ 
 
$$= \frac{  \Gamma^2 ( 1/2-\mu -i\tau) }{ \pi  \Gamma(-2i\tau)\  \Gamma(2(\mu-i\tau))}  \int_0^\infty     e^{-y} y^{2(\mu-i\tau)-1} \int_{1/N}^N  e^{-yt}  G(t)\  dt  dy = g_N(\tau).$$ 
 Thus we find the following representation for the sequence $\{f_N\}_{N\ge 1}$
 
  $$   f_N(\tau) = {2\over \pi} \int_0^\infty      {\rm Re} \bigg[  \frac{  \Gamma^2 ( 1/2-\mu -i\tau)    }{  \Gamma(-2i\tau)  \Gamma(2(\mu-i\tau))}  {}_1F_2 \left( {1\over 2} +\mu -i\tau ;  \ 1-2i\tau, \   \mu-i\tau;\  {y^2\over 4} \right) \ y^{- 2i\tau} \bigg] $$

$$\times e^{-y} y^{2\mu-1} \left(D_{-}^{2\mu} F_{\mu, N} f \right) (y)  dy.\eqno(2.30)$$
 To end this section we show the convergence of the sequence $\{f_N\}_{N\ge 1}$ in the mean square sense to $f(\tau)$ given by the integral
 
 $$     f(\tau) = {2\over \pi} \int_0^\infty      {\rm Re} \bigg[  \frac{  \Gamma^2 ( 1/2-\mu -i\tau)    }{  \Gamma(-2i\tau)  \Gamma(2(\mu-i\tau))}  {}_1F_2 \left( {1\over 2} +\mu -i\tau ;  \ 1-2i\tau, \   \mu-i\tau;\  {y^2\over 4} \right) \ y^{- 2i\tau} \bigg] $$

$$\times e^{-y} y^{2\mu-1} \left(D_{-}^{2\mu} F_{\mu} f \right) (y)  dy\eqno(2.31)$$
 which can be treated as the inversion formula for the  transformation (1.5).  Indeed,  we have, taking into account the above estimates
 
 $$|| f_N- f||_{L_2(\mathbb{R}_+; |\Gamma(2i\tau)/ |\Gamma(1/2-\mu+ i\tau)|^2|^2 d\tau) } \le ||g_N(\tau) -g(\tau)||_{L_2(\mathbb{R}_+; |\Gamma(2i\tau)/ |\Gamma(1/2-\mu+ i\tau)|^2|^2 d\tau) } $$
 
 $$+ ||g_N(-\tau) -g(-\tau)||_{L_2(\mathbb{R}_+; |\Gamma(2i\tau)/ |\Gamma(1/2-\mu+ i\tau)|^2|^2 d\tau) } $$

$$+ {1\over \pi} \left(\int_0^\infty  {1\over |\Gamma(2(\mu-i\tau))|^2}   \bigg| \sum_{n=1}^\infty {(1/2+\mu-i\tau)_n \Gamma( 2(\mu-i\tau+n)) \over 4^n n! \ (1-2i\tau)_n\  (\mu-i\tau)_n} \right.\bigg.$$

$$\left.\bigg. \times  \left( \int_0^{1/N} + \int_N^\infty \right)  {G(t) \ dt \over (t+1)^{2(\mu-i\tau+n)} } \bigg|^2 d\tau\right)^{1/2}$$

 $$+ {1\over \pi} \left(\int_0^\infty  {1\over |\Gamma(2(\mu+i\tau))|^2}   \bigg| \sum_{n=1}^\infty {(1/2+\mu+i\tau)_n \Gamma( 2(\mu+i\tau+n)) \over 4^n n! \ (1+2i\tau)_n\  (\mu+i\tau)_n} \right.\bigg.$$

$$\left.\bigg. \times  \left( \int_0^{1/N} + \int_N^\infty \right)  {G(t) \ dt \over (t+1)^{2(\mu+i\tau+n)} } \bigg|^2 d\tau\right)^{1/2}$$

$$=  ||g_N(\tau) -g(\tau)||_{L_2(\mathbb{R}_+; |\Gamma(2i\tau)/ |\Gamma(1/2-\mu+ i\tau)|^2|^2 d\tau) } $$
 
 $$+ ||g_N(-\tau) -g(-\tau)||_{L_2(\mathbb{R}_+; |\Gamma(2i\tau)/ |\Gamma(1/2-\mu+ i\tau)|^2|^2 d\tau) } $$

$$+ {1\over \pi} \left(\int_0^\infty  \bigg| \bigg[ {}_2F_1\left( {1\over 2}+\mu-i\tau,\  {1\over 2}+\mu-i\tau;\ 1-2i\tau;\ {1\over (t+1)^2} \right) -1\bigg] \bigg.\right.$$

$$\left.\bigg. \times  \left( \int_0^{1/N} + \int_N^\infty \right)  {G(t) \ dt \over (t+1)^{2(\mu-i\tau)} } \bigg|^2 d\tau\right)^{1/2}$$

 $$+ {1\over \pi} \left(\int_0^\infty   \bigg| \bigg[ {}_2F_1\left( {1\over 2}+\mu+i\tau,\  {1\over 2}+\mu+i\tau;\ 1+2i\tau;\ {1\over (t+1)^2} \right) -1\bigg]\bigg.\right.$$

$$\left.\bigg. \times  \left( \int_0^{1/N} + \int_N^\infty \right)  {G(t) \ dt \over (t+1)^{2(\mu+i\tau)} } \bigg|^2 d\tau\right)^{1/2}.$$
 Meanwhile, the  slightly corrected Entry 2.16.33.1 in [1, Vol. II]  implies the representation of the latter Gauss hypergeometric functions
 
 $$ {}_2F_1\left( {1\over 2}+\mu-i\tau,\  {1\over 2}+\mu-i\tau;\ 1-2i\tau;\ {1\over (t+1)^2} \right) = {4^{1-\mu} \Gamma(1+2\mu) \over |\Gamma(1/2+\mu+i\tau)|^4}$$
 
 $$\times  \left(1- {1\over (t+1)^2} \right)^{1/2+\mu+i\tau}  \int_0^\infty x^{2\mu} K_0\left( x \sqrt{ 1- {1\over (t+1)^2} } \right) K_{2i\tau} (x) dx.$$
Hence, employing again the Cauchy-Schwarz  inequality and Entry 2.16.33.2 in [1, Vol. II], we derive

$$\left|  {}_2F_1\left( {1\over 2}+\mu-i\tau,\  {1\over 2}+\mu-i\tau;\ 1-2i\tau;\ {1\over (t+1)^2} \right)\right| \le  {4^{1-\mu} \Gamma(1+2\mu) \over |\Gamma(1/2+\mu+i\tau)|^4} \left(1- {1\over (t+1)^2} \right)^{1/2+\mu} $$

$$\times \left( \int_0^\infty x^{2\mu}\  K^2_0\left( x \sqrt{ 1- {1\over (t+1)^2} } \right) dx\right)^{1/2} \left( \int_0^\infty x^{2\mu}  K^2_{2i\tau} (x) dx\right)^{1/2} $$

$$=   { \Gamma^{2} (1/2+\mu) \   |\Gamma(1/2+\mu+ 2i\tau)|^2 \over |\Gamma(1/2+\mu+i\tau)|^4} \left(1- {1\over (t+1)^2} \right)^{1/4+\mu/2},$$
and for big $\tau$ via Stirling's asymptotic formula 

$$ {|\Gamma(1/2+\mu+ 2i\tau)|^2 \over |\Gamma(1/2+\mu+i\tau)|^4} = O(\tau^{-2\mu}),\quad \tau \to \infty.$$
Therefore for $\mu \in [0, 1/4)$ it has 

$$|| f_N- f||_{L_2(\mathbb{R}_+; |\Gamma(2i\tau)/ |\Gamma(1/2-\mu+ i\tau)|^2|^2 d\tau) } \le 2 ||g_N(\tau) -g(\tau)||_{L_2(\mathbb{R}_+; |\Gamma(2i\tau)/ |\Gamma(1/2-\mu+ i\tau)|^2|^2 d\tau) } $$
 
 $$+ 2 ||g_N(-\tau) -g(-\tau)||_{L_2(\mathbb{R}_+; |\Gamma(2i\tau)/ |\Gamma(1/2-\mu+ i\tau)|^2|^2 d\tau) } $$

$$+ \Gamma^{2} (1/2+\mu)\ { 2^{1/2+\mu} \over \pi}  \  \left(\int_0^\infty\bigg| {\Gamma^2 (1/2+\mu+ 2i\tau) \over \Gamma^4(1/2+\mu+i\tau)} \bigg|^2 \right.$$

$$\left.\bigg. \times  \left( \int_0^{1/N} + \int_N^\infty \right)  { G(t) \ dt \over (t+1)^{2(\mu- i\tau)} } \bigg|^2 d\tau\right)^{1/2}$$

 $$+  \Gamma^{2} (1/2+\mu) \ { 2^{1/2+\mu} \over \pi} \  \left(\int_0^\infty   \bigg| {\Gamma^2 (1/2+\mu+ 2i\tau) \over \Gamma^4(1/2+\mu+i\tau)} \bigg|^2 \right.$$

$$\left.\bigg. \times  \left( \int_0^{1/N} + \int_N^\infty \right)  {\ G(t) \ dt \over (t+1)^{2(\mu+ i\tau)} } \bigg|^2 d\tau\right)^{1/2}$$

$$\le D\bigg[  ||g_N(\tau) -g(\tau)||_{L_2(\mathbb{R}_+; |\Gamma(2i\tau)/ |\Gamma(1/2-\mu+ i\tau)|^2|^2 d\tau) } \bigg.$$
 
 $$+ \bigg.  ||g_N(-\tau) -g(-\tau)||_{L_2(\mathbb{R}_+; |\Gamma(2i\tau)/ |\Gamma(1/2-\mu+ i\tau)|^2|^2 d\tau) } \bigg],\quad \to 0,\quad N \to \infty,$$
where $D >0$ is a constant. Thus we proved 

{\bf Theorem 3}. {\it Let $\mu \in [0, 1/4)$. Then under conditions of Theorem $2$ the inversion formula $(2.31)$ for the transform $(1.5)$ takes place, where $f(\tau) = \lim_{N\to \infty} f_N(\tau)$ with respect to the norm in $L_2(\mathbb{R}_+; |\Gamma(2i\tau)/ |\Gamma(1/2-\mu+ i\tau)|^2|^2 d\tau) $ and the sequence $\{f_N\}_{N\ge 1}$ is defined by $(2.30)$}.

\section{ Special case $\mu=0$}

This section concerns the integral transformation (1.6) of the Lebedev type [5], involving the square of the Macdonald function.   Theorem 1 reads as follows

{\bf Theorem 4}. {\it Let $x >0,\  \ \delta \in \left[0,  \pi/2\right) $ and $f \in L_1\left( \mathbb{R}_+; e^{(\pi-2\delta) \tau}  d\tau\right)$. Then the transform $(1.6)$ is a bounded map $F:  L_1\left( \mathbb{R}_+; \ e^{(\pi-2\delta) \tau}  d\tau\right) \to C\left([x_0, \infty)\right), x_0 > 0$ and can be represented as a composition of the Laplace and Olevskii transforms, namely,
 
$$(F_0 f) (x)=    \int_0^\infty    e^{- x t} G_0(t) dt,\eqno(3.1)$$
where}

$$G_0(t)=   \int_0^\infty  \ {}_2F_1 \left( {1\over 2} -i\tau,\  {1\over 2} + i\tau ; 1 ; - t^2- 2t \right) f(\tau) d\tau.\eqno(3.2)$$
The Olevskii transform (3.2), in turn,  is a bijective map from $L_2(\mathbb{R}_+; (\coth(\pi\tau)/\tau) d\tau) $ onto $L_2(\mathbb{R}_+;  (t+1) dt)$, where  the inverse operator is given by the formula
 $$ f(\tau)= 4\tau \tanh(\pi\tau)  \int_0^\infty  \ (t+1)\ {}_2F_1 \left( {1\over 2} -i\tau,\  {1\over 2} + i\tau ; 1 ;  - t^2-2t \right) G_0\left( t\right) dt,\eqno(3.3)$$
where the convergence is with respect to the norm in $L_2(\mathbb{R}_+; (\coth(\pi\tau)/\tau) d\tau) $.  Moreover, the Parseval equality holds 
$$\int_0^\infty  \   |f(\tau)|^2  \coth(\pi\tau)\  {d\tau\over \tau}  = {4\over \pi} \int_0^\infty \left| G\left(t \right) \right|^2  (t+1) \ dt.\eqno(3.4)$$
The inversion formula for the transformation (1.6) is given by the corresponding analog of (2.31), letting $\mu=0$.   Therefore we have, making a simple substitution  

$$     f(\tau) = 8 \int_0^\infty      {\rm Re} \bigg[ {1\over \Gamma^2(-i\tau)  }\   {}_1F_2 \left( {1\over 2}  -i\tau ;  \ 1-2i\tau, \   -i\tau;\  y^2 \right) \ \left({y^2\over 4}\right)^{- i\tau} \bigg] $$

$$\times e^{-2y} \  (F_{0} f)  (2y)  {dy\over y}\eqno(3.5)$$
with the convergence of the integral (3.5) with respect to the norm in $L_2(\mathbb{R}_+; (\coth(\pi\tau)/\tau) d\tau) $. Fortunately, the kernel in (3.5) can be evaluated in terms of the modified Bessel functions by virtue of the Entry  7.14.1.4 in [1, Vol. III]. Consequently,  it yields

$$ {}_1F_2 \left( {1\over 2}  -i\tau ;  \ 1-2i\tau, \   -i\tau;\  y^2 \right) = y\  \Gamma (-i\tau)  \Gamma (1-i\tau) \left({y^2\over 4}\right)^{i\tau} I_{-i\tau} \left(y\right) I_{-i\tau-1} \left(y\right)$$

$$-  \Gamma^2 (1-i\tau) \left({y^2\over 4}\right)^{i\tau} I^2_{-i\tau} \left(y\right)$$
and correspondingly,

$$ {1\over \Gamma^2(-i\tau)  }\   {}_1F_2 \left( {1\over 2}  -i\tau ;  \ 1-2i\tau, \   -i\tau;\  y^2\right) \ \left({y^2\over 4}\right)^{- i\tau} $$

$$ =  - i\tau \ I_{-i\tau} \left(y\right)  \bigg[  y\  I_{-i\tau-1} \left(y\right)  + i \tau   I_{-i\tau} \left(y\right) \bigg].$$
Employing the recurrence relation for the modified Bessel functions [1, Vol. II], the latter formula can be read as follows  

 $$ {1\over \Gamma^2(-i\tau)  }\   {}_1F_2 \left( {1\over 2}  -i\tau ;  \ 1-2i\tau, \   -i\tau;\  y^2 \right) \ \left({y^2\over 4}\right)^{- i\tau} $$

$$ =  - i\tau y \ I_{-i\tau} \left(y\right) \ {d\over dy} \left[ I_{-i\tau} \left(y\right)\right] = -  i\tau  {y\over 2}  \  {d\over dy} \left[ I^2_{-i\tau} \left(y\right)\right].$$
 Thus the inversion formula (3.5) takes a final form 
 
 $$     f(\tau) = - 2 i  \tau \int_0^\infty   {d\over dy}  \bigg[ I^2_{-i\tau} \left(y\right) - I^2_{i\tau} \left(y\right) \bigg] e^{-2y} \  (F_{0} f)  (2y) \ dy.\eqno(3.6)$$

\bigskip
\centerline{{\bf Acknowledgments}}
\bigskip

\noindent The work was partially supported by CMUP, which is financed by national funds through FCT (Portugal)  under the project with reference UIDB/00144/2020.

\bigskip
\centerline{{\bf References}}
\bigskip
\baselineskip=12pt
\medskip
\begin{enumerate}

\item[{\bf 1.}\ ] A.P. Prudnikov, Yu.A. Brychkov and O.I. Marichev, {\it Integrals and Series}. Vol. I: {\it Elementary Functions}, Vol. II: {\it Special Functions}, Gordon and Breach, New York and London, 1986, Vol. III : {\it More special functions},  Gordon and Breach, New York and London,  1990.

 \item[{\bf 2.}\ ]  R. Sousa, M.  Guerra and S.  Yakubovich,  {\it Convolution-like structures, differential operators and diffusion processes}. Lecture Notes in Mathematics, 2315. Springer, Cham, 2022.

\item[{\bf 3.}\ ]  J. Wimp,  A class of integral transforms,  {\it  Proc. Edinb. Math. Soc.} {\bf 14} (1964), N 2,  33-40. 

\item[{\bf 4.}\ ]  S. Yakubovich and Yu. Luchko, The Hypergeometric Approach to Integral Transforms and Convolutions, {\it Kluwer
Academic Publishers, Mathematics and Applications.} Vol.287, 1994. 

\item[{\bf 5.}\ ] S. Yakubovich, {\it Index Transforms}, World Scientific Publishing Company, Singapore, New Jersey, London and
Hong Kong, 1996.

\item[{\bf 6.}\ ]  S. Yakubovich,  On the Plancherel theorem for the Olevskii transform.  {\it Acta Math Vietnam}.  Vol. {\bf 31} (2006),  249-260.

\end{enumerate}

\vspace{5mm}

\noindent S.Yakubovich\\
Department of  Mathematics,\\
Faculty of Sciences,\\
University of Porto,\\
Campo Alegre st., 687\\
4169-007 Porto\\
Portugal\\
E-Mail: syakubov@fc.up.pt\\

\end{document}